\documentclass[11pt]{amsart}
\usepackage{amssymb}
\usepackage[dvips]{graphicx}

\newcommand{\SD}{{\mathcal{D}}}

\newcommand{\SP}{{\mathcal{P}}}

\newcommand{\Z}{\mathbb{Z}}
\newcommand{\C}{\mathbb{C}}

\newcommand{\R}{\mathbb{R}}

\newtheorem{proposition}{Proposition}[section]
\newtheorem{theorem}[proposition]{Theorem}
\newtheorem{definition}[proposition]{Definition}
\newtheorem{lemma}[proposition]{Lemma}
\newtheorem{conjecture}[proposition]{Conjecture}
\newtheorem{corollary}[proposition]{Corollary}
\newtheorem{remark}[proposition]{Remark}

\begin{document}

\title
{A class of non-fillable contact structures}

\subjclass{Primary: 57R17 Secondary: 53D10.}
\date{October, 2006.}
\keywords{contact structures, fillings}

\author{Francisco Presas}
\address{Departamento de Matem\'aticas \\
Consejo Superior de Investigaciones Cient\'{\i}ficas
\\ 28006 Madrid \\ Spain}
\email{fpresas@imaff.cfmac.csic.es}


\renewcommand{\theenumi}{\roman{enumi}}

\begin{abstract}
A geometric obstruction, the so called "plastikstufe", for a contact structure to not being fillable has been found in \cite{Kl06}. This generalizes somehow the concept of overtwisted structure to dimensions higher than $3$. This paper elaborates on the theory showing a big number of closed contact manifolds having a "plastikstufe". So, they are the first examples of non-fillable high dimensional closed contact manifolds. In particular we show that $S^3 \times \prod_j \Sigma_{j}$, for $g(\Sigma_j)\geq 2$, possesses this kind of contact structure and so any connected sum with those manifolds also does it.
\end{abstract}

\maketitle

\section{Introduction} \label{sect:introduction}
Since the mid-eighties there have been two clear categories of $3$-dimensional contact manifolds: overtwisted and tight. The first ones were introduced by Eliashberg, following Gromov ideas, and happened to satisfy a kind of $h$-principle making their study a mere homotopical question \cite{El89}. Tightness has been more evasive being a coarse classification a matter of only very recent research and being the study purely "geometrical" with no hope of $h$-principle flexibility.

A key result to understand this overtwisted class is a corollary of Gromov's pseudoholomorphic curves foundational paper \cite{Gr85}:
\begin{theorem} \label{thm:gro}
Every overtwisted manifold is non-fillable.
\end{theorem}
Several definitions can be given for fillability. The most restrictive one is Stein-fillable. A contact manifold is Stein-fillable if it can be expressed as a level hypersurface for a proper plurisubharmonic function on a Stein manifold. It is very easy to check that locally, forgetting about the complex structure, a neighborhood of the manifold is symplectomorphic to a neighborhood of the $1$-section of the symplectization $(C\times (0, \infty), d(e^t\alpha))$. This provides another definition. A contact manifold $(C, \alpha)$ is symplectically fillable if there exists a symplectic manifold $(M, \omega)$ with boundary $\partial M$ such that a neighborhood $U$ of the boundary is symplectomorphic to $(C\times (1-\epsilon, 1], d(e^t\alpha))$. There is also a concept of weak filling where the contact manifold $(C, \alpha)$ is the boundary of a symplectic manifold $(M, \omega)$ satisfying that $\ker \alpha$ is a symplectic subspace for $\omega$. This is usually referred to as $\alpha$ being dominated by $\omega$. Again this is a weaker condition than being symplectically fillable. Also, in dimension $3$, it has been proved that these inclusions are strict, so we have
\begin{equation*}
  \text{Stein-fillable} \varsubsetneq \text{Symplectically-fillable} \varsubsetneq \text{Weakly-fillable} \varsubsetneq \text{Tight}.
\end{equation*}
Theorem \ref{thm:gro} works for weakly symplectic fillings. Finally the division \emph{tight-overtwisted} does not translate in a completely faithful way to the fillings since some examples of tight manifolds not having fillings at all have been found . However, in most of the cases tight manifolds are fillable.

All the previous theory has been established only in dimension $3$, since many of the examples of contact structures in higher dimension are expected to be Stein-fillable. The main reason is that a good generalization of the concept of overtwisted to higher dimensions has been hard to find in the last 20 years. The first successful attempt has been \cite{Kl06}. There, it has been shown that the existence of a special geometrical structure on a contact manifold automatically implies the non-fillability of it. This generalizes Theorem \ref{thm:gro} and shows that the definition makes some sense. Let us show it.
\begin{definition} \label{def:plastik}
A plastikstufe in a contact manifold $(C^{2n-1}, \alpha)$ is an embedding $\SP(S)=D^2 \times S \hookrightarrow C$, where $S$ is a closed $(n-2)$-dimensional manifold satisfying that
\begin{enumerate}
\item $B=\{0 \} \times S$ is tangent to the contact distribution $D= \ker \alpha$.
\item $S_D=T\SP(S)\bigcap D$ defines a singular distribution on $\SP(S)$ which has a singular set defined by $B$.
\item $S_D$ is an integrable distribution on $\SP(S)$, in fact a Legendrian foliation.
\item The boundary $\partial \SP(S)= S^1 \times S$ is a closed leaf of the foliation $S_D$, and so it is a Legendrian submanifold.
\item $B$ is a family of singularities of elliptic type, i.e any small $2$-dimensional disk $D_{\epsilon}$, transverse to $B$, has an induced distribution $D\bigcap TD_{\epsilon}$ with a singularity of elliptic type.
\item The rest of the leaves of the foliation $S_D$ are of the form $S \times (0,1)$ satisfying that their boundary is contained in $B\bigcup \partial \SP(S)$
\end{enumerate}
A manifold with a plastikstufe is called overtwisted.
\end{definition}
\begin{figure}[htp]
\includegraphics[scale=0.55]{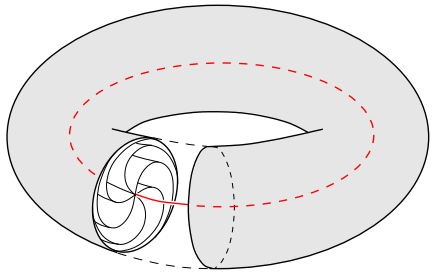} \\
\caption{{\small $\SP(S^1)$ plastikstufe in a $5$-dimensional manifold.}} \label{fig:plastik}
\end{figure}
\vspace{12mm}
Recall that a $2$-dimensional plastikstufe is just an overtwisted disk. So far, it has been proved:
\begin{theorem}[Theorem 1 in \cite{Kl06}] \label{thm:plas_non}
An overtwisted manifold does not admit a semi-positive symplectic filling.
\end{theorem}
We mean by a semi-positive filling a weak symplectic filling given by a semi-positive symplectic manifold. This covers in particular Stein fillings and exact symplectic ones. It is expected to have a result covering all the possible fillings. The main obstruction is technical and is due to not have used "virtual cycles" in the proof of Theorem \ref{thm:plas_non}. In the particular case of dimension $5$ all the fillings can be made "semi-positive" and so the Theorem works for any filling. Remark that in \cite{Gr85} some sketches of the proof of a result similar to Theorem \ref{thm:plas_non} were already given. Moreover, Yuri Chekanov also has an alternative proof to Theorem \ref{thm:plas_non}.

This paper has the task of showing the existence of a very ample class of overtwisted contact manifolds. The main result we will prove is
\begin{theorem} \label{thm:main}
Assume that a contact manifold $(M', \alpha)$ admits a codimension $2$ contact submanifold $N$, with trivial normal bundle, such that the restricted contact form is overtwisted, then $M' \bigcup_N (N \times T^2)$ admits an overtwisted contact structure.
\end{theorem}
Later we will see that using Theorem \ref{thm:main} is simple to construct quite a lot families of high dimensional examples. For instance
\begin{corollary} \label{coro:main}
The manifold $S^3 \times \prod_i \Sigma_i$, where $\Sigma_i$ is a Riemann surface with genus at least $2$, admits an overtwisted contact structure.
\end{corollary}
And so we immediately obtain
\begin{corollary} \label{coro:main2}
The connected sum of any contact manifold with $S^3 \times \prod_i \Sigma_i$ admits an overtwisted contact structure.
\end{corollary}
In dimension $5$, we obtain
\begin{corollary} \label{coro:main3}
The manifold $M_o=(S^3 \times T^2) \bigcup_{S^3} S^5$ admits an overtwisted contact structure.
\end{corollary}
Recall that this manifold has homology $H_2(M_o)=H_3(M_o)=0$ and $H_1(M_o)=H_4(M_o)= \Z^2$. So, it is not a sphere but it has a simple homology type.
It is reasonable to conjecture:
\begin{conjecture}
$S^{2n+1}$ admits an overtwisted structure and so any contact manifold admits an exotic overtwisted structure.
\end{conjecture}
This would partially recover the familiar $3$-dimensional picture in which for any contact structure, it is possible to construct an overtwisted one homotopically equivalent as a plane distribution. Recall that Corollary \ref{coro:main3} gets close to this conjecture in dimension $5$.

Finally observe that no examples of non-fillable contact structures had been found in dimension higher than $3$. All the previous results are a huge source of such structures.

In Section \ref{sect:contact_conn} we define a notion of parallel transport for contact fibrations which is a key step in our proof. The relation between a contact fibration and the group of Hamiltonian contactomorphisms is important for our purposes and it is studied in Section \ref{sect:contact_ham}. Next, we review some old and new constructions of contact structures in Section \ref{sect:constructions}. Section \ref{sect:proof} just places together all the previous results to prove the main theorem. Finally some Corollaries of the main theorem are obtained in Section \ref{sect:remarks}. In particular we give examples of high dimensional manifolds which admit contact structures both fillable and non-fillable lying in the same homotopy class of hyperplane fields.

{\bf Acknowledgments:} \newline
Special thanks to K. Niederkr\"uger whose interest in this project has been amazing. His corrections have made this paper change a lot. Figure \ref{fig:plastik} is due to him. Also I want to thank Maica Mata for her help and encouragement all over these years. Finally I want to thank Vicente Mu\~noz for all the discussions we have had in the last few years.

\section{Contact connections.} \label{sect:contact_conn}
\subsection{Review of the symplectic case.}
We briefly review a classical construction in symplectic geometry as a motivation for our proofs.
\begin{definition}
Let $\pi: M \to B$ be a fibration such that $M$ admits a symplectic structure $\omega$ satisfying that the fibers $F=\pi^{-1}(b)$ of $\pi$ are symplectic submanifolds, then the distribution $H=(TF)^{\bot \omega}$ is called the
symplectic connection associated to the symplectic fibration $(M, B, \pi, \omega).$
\end{definition}
The main result is
\begin{proposition} \label{propo:symp_mono}
Let $\pi: M\to B$ be a symplectic fibration with associated symplectic connection $H$. For each $p_0\in M$, $H$ allows us to lift any path $\gamma:[0,1] \to B$, such that $\pi(p_0)=\gamma(0)$, to a unique path $\tilde{\gamma}_{p_0}:[0,1] \to M$ satisfying:
\begin{itemize}
\item $\tilde{\gamma}_{p_0}(0)=p_0$,
\item $\pi\circ \tilde{\gamma}_{p_0}= \gamma$,
\item $\tilde{\gamma}'_{p_0}(t) \in H_{\tilde{\gamma}_{p_0}(t)}$.
\end{itemize}
Moreover, for $\gamma$ immersed, the induced map
\begin{eqnarray*}
m: F_0=\pi^{-1}(\gamma(0)) & \to & F_1=\pi^{-1}(\gamma(1)) \\
p & \to & \tilde{\gamma}_{p}(1)
\end{eqnarray*}
is a symplectomorphism.
\end{proposition}
\noindent {\bf Proof.} \newline
The first part of the statement, as usual for any connection on a fibration, consists of solving the equation:
\begin{eqnarray*}
& & \pi_*(\tilde\gamma'_{p_0}(t))= \gamma'(t), \\
& & \tilde{\gamma}'_{p_0}(t) \in H.
\end{eqnarray*}
From $\tilde{\gamma}'$ we obtain $\tilde{\gamma}$ as a consequence of the theorem of existence and uniqueness of solutions of ODEs.

For the second part we define the submanifold with boundary $P=\pi^{-1}(\gamma([0,1]))$ that is an immersed submanifold whenever $\gamma$ is immersed. The vector field $X=\gamma_*(\frac{d}{dt})$ is a vector field on $P$ that preserves the fibers $F_t=\pi^{-1}(\gamma(t))$. We want to show that this field infinitesimally preserves the symplectic form, this is
\begin{equation} \label{eq:preserv_symp}
L_X \omega_{|P} = 0.
\end{equation}
Using the Cartan formula we obtain
$$
L_X \omega_{|P} = d i_X \omega_{|P} + i_X d\omega_{|P} = d i_X \omega_{|P}.
$$
Now recall that $ X\in H$ and $H$ is symplectically orthogonal to $TF_t$, so $i_X \omega_{|P}(v)=0$, $\forall v\in TF_t$. Finally we have $i_X \omega(X)=0$ and therefore $d i_X \omega_{|P}=0$ as we wanted to show. So, being true formula (\ref{eq:preserv_symp}), we obtain that $m$ is a symplectomorphism. \hfill $\Box$

\subsection{The contact case.}
We want to adapt the previous language to the contact situation. We denote $D= \ker \alpha$.
\begin{definition}
Let $\pi: C \to B$ be a fibration such that $C$ admits a contact structure $\alpha$ satisfying that the fibers $F=\pi^{-1}(b)$ of $\pi$ are contact submanifolds, then the distribution $H=(TF\bigcap D)^{\bot d\alpha}\subset D$ is called the contact connection associated to the contact fibration $(C, B, \pi, \alpha)$.
\end{definition}
Recall that $D$ is transverse to the fibers if they are contact and so the previous definition makes sense. Again we have that the natural monodromy preserves the contact structure on the fibres:
\begin{proposition} \label{propo:conta_mono}
Let $\pi: C\to B$ be a contact fibration with associated contact connection $H$. For each $p_0\in C$, $H$ allows us to lift any path $\gamma:[0,1] \to B$, such that $\pi(p_0)=\gamma(0)$, to a unique path $\tilde{\gamma}_{p_0}:[0,1] \to M$ satisfying:
\begin{itemize}
\item $\tilde{\gamma}_{p_0}(0)=p_0$,
\item $\pi\circ \tilde{\gamma}_{p_0}= \gamma$,
\item $\tilde{\gamma}'_{p_0}(t) \in H_{\tilde{\gamma}_{p_0}(t)}$.
\end{itemize}Moreover, for $\gamma$ immersed, the induced map
\begin{eqnarray*}
m: F_0=\pi^{-1}(\gamma(0)) & \to & F_1=\pi^{-1}(\gamma(1)) \\
p & \to & \tilde{\gamma}_{p}(1)
\end{eqnarray*}
is a contactomorphism.
\end{proposition}
\noindent {\bf Proof.} \newline
Again the first part of the statement is the usual lift defined for any connection on a fibration.

For the second part we mimic the proof of Proposition \ref{propo:symp_mono} defining the submanifold with boundary $P=\pi^{-1}(\gamma([0,1]))$. It is an immersed submanifold whenever $\gamma$ is immersed. The vector field $X=\gamma_*(\frac{d}{dt})$ is a vector field on $P$ that preserves the fibers $F_t=\pi^{-1}(\gamma(t))$. We want to show that this field infinitesimally preserves the contact structure, that is equivalent to
\begin{equation} \label{eq:preserv_cont}
L_X \alpha_{|P} = f\alpha_{|P},
\end{equation}
for some function $f$. We claim that $f$ has to be chosen to be $f=d\alpha(X, R_{F_t})$, where $R_{F_t}$ is the Reeb vector field associated to the fiber $F_t$. Using the Cartan formula we obtain
$$
L_X \alpha_{|P} = d i_X \alpha_{|P} + i_X d\alpha_{|P} = i_X d \alpha_{|P},
$$
because $X\in D= \ker \alpha$.
Now recall that $ X\in H$ and $H$ is symplectically orthogonal to $TF_t\bigcap D$, so $(i_X d\alpha_{|P})(v)=0$, $\forall v\in TF_t\cap D$. We also have $i_X d\alpha(X)=0$. So we have
that the equation (\ref{eq:preserv_cont}) is true for all $v\in D\cap TP$. To finish we check it against $R_{F_t}$ and we obtain
$$
(L_X \alpha) (R_{F_t})= d\alpha (X, R_{F_t}) = d\alpha (X, R_{F_t}) \alpha (R_{F_t})= f \alpha(R_{F_t}). $$
Therefore the formula (\ref{eq:preserv_cont}) is true and $m$ is a contactomorphism. \hfill $\Box$

\begin{remark} \label{rem:important} \- \newline
\begin{enumerate}
\item The contact connection depends on the choice of contact form and so different contact forms supported on the same contact structure give rise to different monodromies. Monodromies are isotopic inside the contactomorphism group for homotopic contact structures.
\item Let us analyze the most simple case. Choose a contact manifold $(C_0, \alpha_0)$ and define on
$C_0\times D^2$ the standard product structure, which in  polar coordinates is written as
\begin{equation} \label{eq:std_prod}
\alpha = \alpha_0 + r^2 d\theta.
\end{equation}
It is simple to check that the Reeb vector field $R_{0}$ restricts to the Reeb vector field of each fiber. Therefore $f=0$ in formula (\ref{eq:preserv_cont}) and so the parallel transport also preserves the contact form itself (not just its kernel). If we take a path $\gamma:[0, 2\pi] \to D^2$ defined as $\theta \to (r(\theta), \theta)$ we obtain that
$$\gamma'=\frac{\partial r}{\partial \theta} \frac{\partial}{\partial r} + \frac{\partial}{\partial \theta}, $$
and the lift of this field is
$$
X= \frac{\partial r}{\partial \theta} \frac{\partial}{\partial r} + \frac{\partial}{\partial \theta}- r^2 R_0.
$$
Thus, the monodromy at time $\theta$ is generated (in the vertical direction) by the vector
field $-r^2 R_0$. Therefore, after integration, the monodromy is given by the flow of the Reeb vector field in the fiber at time
$$t= -\int_0^{2\pi} r^2(\theta) d\theta= -2Area(\gamma), $$
where $Area(\gamma)$ is the area of the domain bounded by the closed curve $\gamma$.
\end{enumerate}
\end{remark}
\subsection{Parallel transport of plastikstufes.}
Our next task is to understand how a plastikstufe behaves under parallel transport. We just want to check two elementary facts.
\begin{lemma} \label{lem:parallel}
Let $\pi:(C, \alpha) \to B$ be a contact fibration with associated contact connection $H$. Let $N$ be an isotropic submanifold of $F_b$ for some $b\in B$ and $\gamma:[0,1] \to B$ an immersed path such that $\gamma(0)=b$. Denote $m_t$ the parallel transport induced by $\gamma$ between $F_0$ and $F_t$, then $m(N)=\bigsqcup_{t=0}^1 m_t(N)$ is an immersed isotropic submanifold of $C$.
\end{lemma}
\noindent {\bf Proof.} \newline
The new manifold $m(N)$ admits the following tangent decomposition $Tm(N)_p=(m_t)_*(TN) \bigoplus \langle
\tilde{\gamma}'(t) \rangle$.
$m_t$ is a contactomorphism and therefore $(m_t)_*(TN)\subset D$. Moreover, the
parallel transport is defined by trajectories following $H$ that is a subspace of $D=\ker \alpha$ and so $\tilde{\gamma}'(t)\in D$. Thus, $Tm(N)\subset D$ and so $m(N)$ is an isotropic submanifold. \hfill $\Box$

\begin{remark} \- \newline
\begin{enumerate}
\item The new submanifold is embedded if $\gamma$ is embedded or if the contactomorphism $m_t$ separates $N$ from itself in all the crossing fibers.
\item If $N$ is Legendrian and $dim B=2$ then $m(N)$ is Legendrian, just because of the dimensions.
\end{enumerate}
\end{remark}

\begin{lemma} \label{lem:induced}
Let $\pi:(C, \alpha) \to B$ be a contact fibration with contact connection $H$ and assume that $\dim B=2$. Let $\SP(S)$ be a plastikstufe in $F_0$ and $\gamma:S^1 \to B$ a closed immersed (resp. embedded) path starting at $b_0$. Assume that the monodromy $m_1: F_0 \to F_1=F_0$ restricts to the identity on the plastikstufe, then $m(\SP(S))= \bigsqcup m_t(\SP(S))$ is an immersed (resp. embedded) plastikstufe with core $S\times S^1$ in $C$.
\end{lemma}
\noindent {\bf Proof.}
Taking into account that the parallel transport follows $D$, it is easy to check that the parallel transport of the distribution, indeed foliation, $D_{\SP(S)}=D_0 \bigcap T\SP(S)$ generates the distribution $D_{m(\SP(S))}= D \bigcap Tm(\SP(S))$. Recall that given a flow $\phi_t$ and an integrable distribution $\SD$ the object $\phi_t(\SD)$ is still an integrable distribution, moreover $\bigsqcup \phi_t(\SD)$ is still integrable. Since the parallel transport is a flow we have that $D_{m(\SP(S))}$ is a foliation (the distribution remains integrable), still by isotropic leaves because of Lemma \ref{lem:parallel}.

On the other hand the core $S$ of $\SP(S)$ generates by parallel transport the core $S \times S^1$ of the new plastikstufe. Finally the boundary transports to create a new Legendrian boundary, again thanks to Lemma \ref{lem:parallel}. \hfill $\Box$

The previous results show that it is very simple to create immersed plastikstufes. In fact, it can be done in a very local picture.

\begin{proposition}
Let $(C, \alpha)$ be a contact manifold and $N \subset C$ an overtwisted codimension $2$ contact submanifold of it with trivial symplectic normal bundle. Then there exists an immersed plastikstufe of $C$ on an arbitrary small neighborhood of $N$.
\end{proposition}
\noindent {\bf Proof.}
As a consequence of the Gray's stability theorem, there exists a neighborhood $U$ of $N$ and a contactomorphism $\phi:U \to N \times D^2$, satisfying that the target space is equipped with either:
\begin{itemize}
\item the standard contact distribution defined as in formula (\ref{eq:std_prod}),
\item the contact distribution defined as
$$ \alpha = \alpha_0 - r^2 d\theta. $$
\end{itemize}
Assume that we are in the first case, we leave the proof of the second one for the reader. We fix the standard form in this fibration to define a contact connection $H$.

Fix a curve $\gamma: S^1 \to D^2$ satisfying:
\begin{enumerate}
\item $\gamma(0)=0$
\item $\gamma(p)=-\gamma(p+\pi)$. (radial symmetry)
\item $\gamma$ is immersed, with a unique self-intersection point $\gamma(0)=\gamma(\pi)=0\in D^2$.
\end{enumerate}
Because of (ii) and (iii), the domain bounded by $\gamma$ has total area $0$. Therefore the monodromy map, according to
Remark \ref{rem:important}:(iii), is the identity. Now we use Lemma \ref{lem:induced} to construct an immersed plastikstufe $\SP(S\times S^1)$ in $N \times D^2$ out of the one contained in the fiber $N \times \{ 0 \}$.
Using $\phi^*$ we get a plastikstufe in a neighborhood of $N$. Recall that $\gamma$ can be chosen arbitrarily close to $\{ 0 \} \in D^2$, thus the plastikstufe can be constructed arbitrarily close to $N$.

\begin{remark}
The plastikstufe would be embedded if the monodromy separates the initial plastikstufe from itself at the crossing fibers. For instance, assume that we are in the $5$ dimensional case, therefore the fiber is $3$ dimensional and the initial plastikstufe is a classical overtwisted disk. In our previous lemma there is only one crossing point at time $1/2$. Recall that the monodromy at that time is given by the flow of the Reeb vector field at a time $t_0$. This time $t_0$ is given by minus twice the area of the domain bounded by the closed curve $\gamma:[0, 1/2] \to D^2$. So any negative value is possible (up to a limit given by half of the area of $D^2$). Now we choose as $\alpha_0$, among all the possible options, one extending the canonical contact form of the canonical neighborhood of the overtwisted disk given by:
$$
\alpha_{ot}= \cos (r)dz + r \sin(r) d\theta. $$
Recall that the overtwisted disk is given by a slight perturbation of $D(2\pi) \times \{ 0 \}$. The Reeb vector field is
$$ R= \frac{\sin r + r\cos r}{r+ \sin r \cos r}\frac{\partial}{\partial z} + \frac{\sin r}{r + \sin r \cos r} \frac{\partial}{\partial \theta} $$
and it is radially invariant. Moreover its $z$-component is positive at $r=0$ and negative at $r=2\pi$. This implies that the image of the disk through the flow at any time intersects the initial disk along a circle. It is impossible to separate the disk using the monodromy from itself. One can check that any other choice of contact form offers the same result as there is a topological obstruction to separate the disk from itself using the Reeb flow for small times.
\end{remark}
The previous result is not unexpected. Also in dimension $3$, it is simple to build immersed overtwisted disks for any contact $3$-dimensional manifold: those disks exist inside the standard contact ball. The problems arise when one tries to construct an embedded one. To do that with a plastikstufe we need to add some new ideas and to change a bit the topological picture.

\section{Contact Hamiltonians as contact monodromies.} \label{sect:contact_ham}
We are going to show how to produce an arbitrary Hamiltonian contactomorphism by parallel transport on a cleverly chosen contact fibration. Again the ideas of these results come from the symplectic fibrations case. However, we do not review that theory because all our results are clear enough without referring them to the symplectic counterparts.

\subsection{Review of basic definitions and results.}
\begin{definition}
The Hamiltonian vector field $X$ associated to the Hamiltonian function $H:C \to \R$ is the unique vector field defined by the pair of equations:
\begin{eqnarray}
i_X \alpha & = & H,  \label{eq:Ham1} \\
i_X d\alpha & = & (d_RH) \alpha -dH. \label{eq:Ham2}
\end{eqnarray}
\end{definition}
A computation shows that for any Hamiltonian field $L_X \alpha= (d_RH) \alpha$, therefore the field preserves the contact structure. Moreover
\begin{theorem}
Any contactomorphism connected to the identity, in the contactomorphism group, is the flow of a time dependent Hamiltonian vector field.
\end{theorem}
\noindent {\bf Proof.}
If the contactomorphism is connected to the identity, it is generated by a time-dependent vector field $X_t$. Choose $H_t=i_{X_t} \alpha$ and check that the associated Hamiltonian contactomorphism coincides with $X_t$.

\begin{lemma}
A contactomorphism $C^2$-close to the identity (inside the contactomorphisms group $Cont(C)$) can be generated by a $C^1$-small time dependent Hamiltonian function.
\end{lemma}
\noindent {\bf Proof.}
Take $\phi_t:C\to C$ a contactomorphism. By hypothesis, it can be generated by a vector field $\{X_t\}_{t=0}^1$ satisfying
$$ \left( \frac{\partial \phi_t(p)}{\partial t} \right)_{t=t_0} = X_t(\phi_t(p)). $$
Thus a bound in the $C^2$-norm of $\phi_t$ provides a bound in the $C^1$-norm of $X_t$. Therefore
the equation (\ref{eq:Ham1}) provides also bounds for the $C^1$-norm of $H_t$. Being precise we obtain
that $||X_t||_{C^2}\leq \delta$ implies that $||H_t||_{C^1} \leq M\delta$, and $M$ is a constant that
only depends on the contact form $\alpha$. \hfill $\Box$

\subsection{Hamiltonians as monodromies}
The following result clarifies a general geometric picture, though we restrict our proof to the very particular case we will need
\begin{proposition} \label{propo:choose_mono}
Let $\pi: (C, \alpha_0) \times D^2(\delta) \to D^2(\delta)$ be the product contact fibration its the standard contact structure. For each $\alpha_0$ and $\delta>0$, there exists a constant $\epsilon>0$ such that for any Hamiltonian function $\{g_t \}_{t=0}^1: C \to \R $ satisfying:
$$
|g| \leq \epsilon, |dg| \leq \epsilon,
$$
where $t$ is considered as another coordinate for the understanding of the second equation, there is a contact form $\overline{\alpha}$ isotopic to the initial one such that the path (given in polar coordinates)
\begin{eqnarray*}
\gamma: [0,\delta] & \to & D^2(\delta) \\
t & \mapsto &  (\frac{1}{2}(t+ \frac12\delta), 0)
\end{eqnarray*}
has as associated monodromy the contactomorphism generated by the Hamiltonian $\{g_t\}$. Moreover the deformation of the contact form can be supported in the set $[\delta/4, 3\delta/4] \times [-\pi/4, \pi/4]$ (this neighborhood is written in polar coordinates).
\end{proposition}
\noindent {\bf Proof.} \newline
Let $\phi_t$ be the contact flow associated to the Hamiltonian $g_t$. Recall that by performing a reparametrization $t=t(s)$ (not necessarily a diffeomorphism) we obtain $\hat{\phi}_s=\phi_{t(s)}$. Just choose the reparametrization in such a way that $\hat{\phi}_s=\phi_0=id$ for $s\in [0, \beta]$ and $\hat{\phi}_s=\phi_1$ for $s\in [1-\beta,1]$ for a small $\beta>0$. The Hamiltonian generating function changes to $\hat{g_s}= \frac{dt}{ds}g_{t(s)}$. Taking $\beta$ small enough this just increases the derivatives of $\hat{g}_s$ slightly. Therefore we may assume that $g_t$ satisfies that $g_t=0$ for $t\in [0, \beta] \cup [1-\beta, 1]$ at the cost of decreasing slightly the constant $\epsilon$ in our statement. This allows us to extend the definition of $g_t$ as $\{g_t \}_{t= -\infty}^{\infty}: C \to \R$ with $g_t=0$ for $t \not\in (\beta, 1-\beta)$.

Let us start with the case $\delta=2$, we solve the general case later. Choose a cut-off function $\xi:\R \to \R$ satisfying:
\begin{enumerate}
\item $\xi(t)=1$ for $t\in(-\pi/8, \pi/8)$.
\item $\xi(t)=0$ for $|t|\geq \pi/4.$
\item $0 \leq \xi(t) \leq 1$, for all $t\in \R$.
\item $|\xi'(t)| \leq 4. $
\end{enumerate}
Define $\tilde{g}(r,\theta,p)=\xi(\theta) \cdot g_{r-1/2}(p)$ which is globally defined in $C\times D^2(2)$. Now we can define
$$\overline{\alpha}= \alpha_0 + r^2 d\theta - \tilde{g}dr. $$
It is obvious that, for $\epsilon$ small enough, this is a contact form isotopic to
the standard contact form $\alpha$, just by lineal interpolation.

The exterior differential of $\overline{\alpha}$ is
$$ d\overline{\alpha}= d\alpha_0 +2rdr\wedge d\theta -d\tilde{g}\wedge dr. $$
We impose the conditions of Proposition \ref{propo:conta_mono} to lift the vector $\frac{\partial}{\partial r}$. This implies that $X=\frac{\partial}{\partial r}+ v$, where $v\in TC$. Moreover we also have:
\begin{itemize}
\item $i_X \overline{\alpha} = 0$, which translates to
\begin{equation}
i_v \alpha_0 = \tilde{g}. \label{eq:first_cond}
\end{equation}
\item $i_X d\overline{\alpha}(w) =0$, for all $w\in D_0$, which translates to
\begin{equation}
i_v d \alpha_0 (w) + d_w \tilde{g} = 0, \- \forall w\in D_0. \label{eq:second_cond}
\end{equation}
\end{itemize}
We can extend the equation (\ref{eq:second_cond}) from $D_0$ to $TC$ obtaining
\begin{equation}
i_v d\alpha_0 = s\alpha -d\tilde{g}, \label{eq:third_cond}
\end{equation}
where $d\tilde{g}$ is the exterior differential of $\tilde{g}$ restricted to $TC$. We substitute the Reeb vector field on equation (\ref{eq:third_cond}) to obtain $s=d_{R_0} \tilde{g}$. Thus, equations (\ref{eq:first_cond}) and (\ref{eq:third_cond}) are the equations for $v$ to be the Hamiltonian vector field associated to the Hamiltonian function $\tilde{g}(r,\cdot, \cdot)$ for $r\in [1/2, 3/2]$. We are restricted to the path, in polar coordinates $[1/2, 3/2] \times \{ 0 \}$ and therefore we have that $\tilde{g}(r, 0, \cdot)=g_{r-1/2}$ that is what we wanted.

Now we go to the case of a general $\delta>0$. We can perform a change of parameter in the associated contact flow $\phi_t$ to $\hat{\phi}_s=\phi_{s/\delta}$, for this we change the generating Hamiltonians to $\hat{g}_s= (g_{s/\delta})/\delta$. We can follow the previous proof from that point since the new Hamiltonian functions $\hat{g}_s$ are now supported for $s\in[0, \delta]$. At the end of the day we obtain that the new contact form $\bar{\alpha}$ generates the required contactomorphism $\phi_1$ as the monodromy generated over the path $[1/2 \delta, 3/2 \delta]\times \{ 0 \}$. The only remark to be made is that the Hamiltonians $\hat{g}$ have as $C^1$-norm the $C^1$-norm of $g$ multiplied by $\delta^{-1}$. Therefore the constant $\epsilon>0$ to be chosen to make $\bar{\alpha}$ isotopic to $\alpha$ depends on $\delta$ as announced in the statement of the Proposition. \hfill $\Box$

\section{Review of some constructions.} \label{sect:constructions}
We need to review three constructions in order to complete our insight of the problem.
\subsection{Contact structures on $C\times T^2$.}
We briefly review the idea of the proof of the following result
\begin{theorem}[Theorem 1 in \cite{Bo03} \label{thm:frederic}]
Let $(C, \alpha)$ a contact manifold. There is a contact structure on $C \times T^2$ which is $T^2$-invariant.
\end{theorem}
\noindent {\bf Proof.}
F. Bourgeois just takes an open book decomposition compatible with $(C, \alpha)$ as constructed in \cite{Gi02}. This object is defined as
\begin{definition} \label{def:openbook}
An open book associated to a manifold $M$ is a pair $(\phi, N)$ such that:
\begin{itemize}
\item $N$ is a codimension $2$ submanifold with trivial normal bundle.
\item $\phi:M-N \to S^1$ is a regular fibration.
\end{itemize}
\end{definition}
\begin{definition} \label{def:compbook}
A compatible open book $(\phi, N)$ of a contact manifold $(C, \alpha)$ is a topological open book satisfying
\begin{itemize}
\item The fibers of $\phi$ are symplectic with respect to the restriction of the closed $2$-form $d\alpha$.
\item The submanifold $N$ is contact.
\item The Reeb vector field is tangent to $N$.
\end{itemize}
\end{definition}
\begin{theorem}[Theorem in \cite{Gi02}] \label{thm_giroux}
Any cooriented contact manifold $(C,D)$ supports a compatible open book for some contact form $\alpha$ associated to $D$.
\end{theorem}
Assume that $\dim C= 2n-3$. To prove Theorem \ref{thm:frederic} we take a compatible open book $(\phi, N)$ and $\phi: N \to S^1 \subset \C$ is multiplied by a function
$\rho$ constructed as a distance to $N$ (constant outside a neighborhood of $N$). So we get a new map $\Phi= \rho \cdot \phi= (\Phi_1, \Phi_2)$. Then the contact form in the product is just
\begin{equation}
\alpha_{\epsilon} = \epsilon(\Phi_1 d\theta_1 + \Phi_2 d\theta_2) + \beta, \label{eq:frede_e}
\end{equation}
where $(\theta_1, \theta_2)$ are lineal coordinates of the torus and $\beta$ is the contact form on $C$. It is a simple computation to write down
\begin{eqnarray*}
\alpha_{\epsilon} \wedge (d \alpha_{\epsilon})^{n-1} & = & -\epsilon^2(n-1)(d\beta)^{n-2} \wedge \rho^2 d\phi \wedge d \theta_1 \wedge d\theta_2 -
\\
& & - \epsilon^2(n-1)(n-2) \beta \wedge (d\beta)^{n-3} \wedge \rho d\rho \wedge d \phi \wedge d\theta_1 \wedge d \theta_2.
\end{eqnarray*}
The two expressions are equal or greater than zero and one is strictly bigger than zero away from $N$ and the other one in a neighborhood of it. In the original paper $\epsilon=1$, but realize that for any $\epsilon \neq 0$ positive or negative the statement is still true.
This will be important for our purposes. \hfill $\Box$

We can easily check that
\begin{lemma} \label{lemma:Frede_bound}
Let $T^2=\R^2 /\Z^2$ be the standard $2$-torus. There is a constant $M>0$ such that for any path $\gamma: [0,\delta] \to T^2$, parametrized by the arc-length satisfying that $|\gamma''|<1$, the monodromy generated by $\gamma$, for the connection associated to $\alpha_{\epsilon}$, after identifying $F_0 \cong C$ and $F_1 \cong C$, is a contactomorphism generated by a Hamiltonian function $(H_t)_{t=0}^1$ satisfying that
$$ |H_t| \leq M |\epsilon| \delta, |dH_t| \leq M |\epsilon| \delta. $$
\end{lemma}
\noindent {\bf Proof.} \newline
It is a straightforward computation. Just take into account that the constant $M$ depends on the norms of $\alpha$, $\Phi_1$, $\Phi_2$ and some of their derivatives. \hfill $\Box$

\subsection{Connected sum of contact fibrations.}
We briefly recall a result of Geiges \cite{Ge97}. We keep the proof, but we adapt the details to our language. The result is
\begin{theorem} \label{thm:geiges}
The fibered sum of contact fibrations over surfaces, satisfying that
\begin{itemize}
\item The fibers of each fibration are contactomorphic,
\item The normal bundles have opposite orientations,
\end{itemize}
is a contact fibration.
\end{theorem}
{\bf Proof.}
Denote $C_j \to B_j$ the contact fibrations with contact forms $\alpha_j$ ($j=1,2$). Denote by $F_0^j= \pi_j^{-1}(b_j)$ the pair of fibers we are going to identify. Moreover, take the connection $H_j$ on the fiber $F_0^j$ as a geometric realization of the normal bundle. We may assume that the fibration can be trivialized to the standard product picture given by
\begin{equation}
\alpha_j= \alpha_0 + \rho(-1)^j r^2 d\theta, \label{eq:normal_model}
\end{equation}
where $\alpha_0$ is the contact form on the fiber and $(r, \theta)$ are polar coordinates on the disk $D^2_j(0, 1)$. This is always possible after small perturbation of the contact form, provided the constant $\rho>0$ is chosen to be small enough (recall that it is usual to assume the radius small enough, but this is equivalent after scaling to fix the radius to $1$ and to introduce $\rho$).
Now we embed $D^2_j(0,1) \subset \R^2 \times \{ (-1)^{j} \} \subset \R^3$. Moreover we construct a hypersurface $S$ on $(\R^+ - \{ 0 \}) \times \R$ (with coordinates $(r,z)$) defined by the zero set of a function $F: (\R^+ - \{ 0 \}) \times \R \to \R$ satisfying that its zero set contains the hyperplanes $\{ z= \pm 1 \}$ for $r\geq 1/2$. Moreover we impose that $F$ is transverse to zero and
\begin{equation}
\left\{ \begin{array}{l} \frac{\partial F}{\partial r} \geq 0 \\
\frac{\partial F}{\partial r} > 0 , z=0, \\
\left( \frac{\partial F}{\partial z} \right) \cdot z \leq 0.
\end{array} \right.
\end{equation}
\begin{figure}[htp]
\includegraphics[scale=0.35]{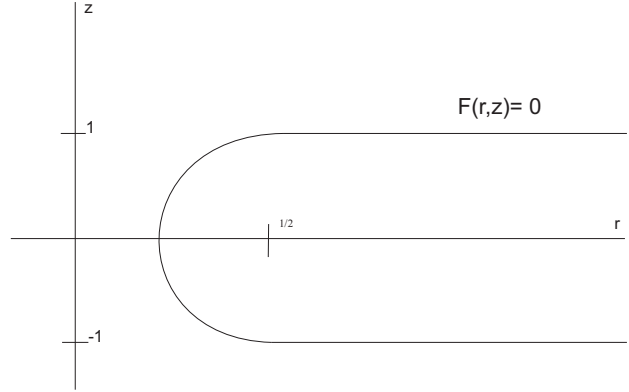} \\
\caption{{\small Hypersurface $S \subset \R^3$}}
\end{figure}
\vspace{12mm}
We define a hypersurface $H\subset \R^3$ by the zero set of a function $G(r,\theta, z) = F(r, z) = 0$, which is obviously rotationally symmetric.
Now the result just states that the form $\hat{\alpha}$ on $F\times \R^3$ given by $\alpha_0 + z \rho r^2 d\theta$ restricted to $F\times H$ is a contact form which restricts to the contact forms $\alpha_j$ on the hypersurfaces $P_{\pm}= \{ z= \pm (-1)^{j+1} \}$. Check \cite{Ge97} for the computation that shows that $\hat{\alpha}$ is contact when restricted to $F\times H$. So, the contact fibration over the annulus $C(\rho)= H \bigcap \{ r \leq 1 \}$ acts as a normal model to construct a contact connection on the connected sum.
Moreover the $C^3$-norm of $G$ can be bounded independently of $\rho$. This implies
that the  $C^2$-norms of the new $\hat{\alpha}$ can be chosen independently of $\rho$. \hfill $\Box$

We can also control the monodromy of the gluing region
\begin{lemma} \label{lemma:Geiges_bound}
The contact connection arising from the connected sum construction of two contact fibrations satisfy that
the holonomy map of a radial path from the interior boundary to the exterior boundary of the annulus $C(\rho)$ can be generated by a Hamiltonian $\{ H_t \}_{t=0}^1$ satisfying $|H_t|_{C^2} \leq M \rho$.
\end{lemma}
{\bf Proof.}
To prove it we just have to check that the contact connection is bounded by a multiple of $\rho$. This is obvious from the construction. \hfill $\Box$

The result really holds for any number of derivatives, but second order bounds are enough for our purposes.

\subsection{An overtwisted $3$-sphere on an exotic $S^5$.}
We are going to define an exotic contact structure on $S^5$ such that it contains a contact $3$-sphere with an overtwisted disk. K. Niederkruger suggested this result to the author and offered a proof, different from the one presented here.
\begin{theorem} \label{propo:exotic_sphere}
There is an exotic contact structure on $S^5$ for which there is a contact sphere $S^3 \subset S^5$ whose contact structure, defined by restriction, is the standard overtwisted one.
\end{theorem}
{\bf Proof.}
We need to introduce some classical results contained in \cite{Mi68}
\begin{theorem} \label{thm:Milnor_fib}
Take $f:U \subset (\C^{n+1}, 0) \to (\C, 0)$ a holomorphic polynomial with an isolated singularity at the origin. Given any $\epsilon > 0$ sufficiently small then:
\begin{enumerate}
\item Define $V=Z(f)$. $K_{\epsilon}= V \bigcap S^{2n+1}(\epsilon)$ is a smooth submanifold of dimension $2n-1$ embedded as a submanifold of the $2n+1$ sphere $S^{2n+1}(\epsilon)$.
\item $\phi_f: \frac{f}{|f|}: S^{2n+1}(\epsilon) - K_{\epsilon} \to S^1$ is a smooth open book.
\end{enumerate}
\end{theorem}
This fibration is called the Milnor fibration of the isolated complex singularity. Moreover there is another canonical fibration
\begin{theorem} \label{thm:Milnor_tube}
Take $f:U \subset (\C^{n+1}, 0) \to (\C, 0)$ a holomorphic polynomial with an isolated singularity at the origin. Given any $\epsilon > 0$ sufficiently small, choose a constant $\delta>0$ small enough with respect to $\epsilon$ such that all the fibers $f^{-1}(t)$ for $|t|\leq \delta$ meet $S^{2n+1}(\epsilon)$ transversally. Let $C_{\delta} = \partial D_{\delta} \subset \C$ the circle of radius $\delta$ centered at $0$, and set $N(\epsilon, \delta)= f^{-1}(C_{\delta}) \bigcap B_{2n+2}(\epsilon)$. Then:
$$ f_{|N(\epsilon, \delta)}: N(\epsilon, \delta) \to C_{\delta} \cong S^1 $$
is a $C^{\infty}$ fiber bundle. Moreover $f_{|N(\epsilon, \delta)}$ and $\phi$ are equivalent fiber bundles.
\end{theorem}
This second fibration is called the Milnor tube. This theorem automatically shows that whenever the singularity is a Morse one, the monodromy is a generalized positive Dehn twist (see \cite{Se99}).

Now we need the inverse of Theorem \ref{thm_giroux} that is much easier to prove. It reads as
\begin{proposition} (\cite{Gi02}) \label{propo:giroux}
Given an open book $(\phi, N)$ satisfying that the leaves are Stein and the monodromy is generated by a symplectomorphism, then there is a contact structure compatible with the open book.
\end{proposition}

Assuming all the previous results, take the function $g(z_1, z_2, z_3)= z_1\bar{z}_2 + z_1z_3 + \bar{z}_2 z_3$. We want to check if it has an associated Milnor fibration, in spite of not being holomorphic. We define
\begin{eqnarray*}
e: \C^3= \C \times \C \times \C & \to & \C \times \C \times \C = \C^3 \\
(z_1, z_2, z_3) & \to & (z_1, \bar{z}_2, z_3).
\end{eqnarray*}
We have that $g \circ e$ is a holomorphic polynomial with a Morse singularity at the origin, therefore Theorems \ref{thm:Milnor_fib} and \ref{thm:Milnor_tube} apply and we obtain that $\phi_{g \circ e}$ is an open book. Moreover recall that the page of the open book $\phi_{g \circ e}(-1)$ is symplectic and the monodromy is a generalized Dehn twist, and so it can be represented by a symplectomorphism. Therefore the open book defines a contact structure in $S^5$, the standard one. Now we want to study $\phi_{g}$. The link (binding) is exactly the same as for $\phi_{g \circ e}$. The page is defined as $e(\phi_{g \circ e}(-1))$. It is diffeomorphic to the page $\phi_{g \circ e}(-1)$ and so it admits a symplectic structure. Theorem \ref{thm:Milnor_tube} allows us to compute the monodromy (just reflecting through $e$) and it is obtained a generalized Dehn twist along a Lagrangian sphere as in the holomorphic case. The only difference is that the orientation of the twist is reversed. Recall that a reverse generalized Dehn twist admits a representation by a symplectomorphism. This proves that $\phi_g$ is an open book that supports a contact structure, thanks to Proposition \ref{propo:giroux} (as we will later see, it is not equivalent to the standard one).

Next we restrict our construction to $\C^2 \subset \C^3$. We have that $f_{|\C^2}= z_1 \bar{z}_2$. It is well known that the associated Milnor fibration for that case is the standard fibration for the complementary of the negative Hopf link \cite{Et05}. This open book produces the standard overtwisted structure in $S^3$. Recall that the pages of this open book can be understood as symplectic submanifolds of the pages of the one generated by $g$. Moreover the generalized Dehn twist can be arranged to preserve the page of $S^3$ when restricted to it. This makes the standard overtwisted $S^3$ a contact submanifold of $S^5$. \hfill $\Box$

\section{Proof of the main results.} \label{sect:proof}

\subsection{Proof of Theorem \ref{thm:main}.}
By Gray's stability theorem there is a neighborhood of $(N, \alpha_{ot})$ contactomorhpic to $(N \times D^2(\delta), \alpha_{ot} \pm r^2 d\theta)$ for a fixed $\delta>0$. We assume a positive signum in the previous expression, being the negative signum left for the careful reader. Now we recall Theorem \ref{thm:frederic} to construct a contact structure $\alpha_+$ on $N \times T^2$, where $N$ is equipped with the overtwisted structure. Moreover in the equation (\ref{eq:frede_e}) defining the contact structure we take $\epsilon$ to be negative and small enough. The precise value of the constant $\epsilon$ will be adjusted later.

Recall Lemma \ref{lem:induced}, this allows us to create an immersed plastikstufe by parallel transport over a path $\gamma:S^1 \to D^2(\delta)$. Let us describe the first half of the path, being the other half defined by symmetry with respect to the origin. The path is defined in $D^2(\delta) \subset \R^2$ as the union of
\begin{enumerate}
\item $[0, 3\delta/4] \times \{ 0 \}$,
\item  A path, in polar coordinates
\begin{eqnarray*}
\gamma: [0, \pi/2] & \to & \R^+ \times [0, 2 \pi] \\
\theta & \mapsto & (r(\theta), \theta),
\end{eqnarray*}
satisfying
\begin{enumerate}
\item $r(0)= r(\pi/2) = 3\delta/4$,
\item $3\delta/4 <  r(\theta) < \delta$, for $\delta \in (0, \pi/2)$,
\end{enumerate}
\item $\{ 0 \} \times [0, 3\delta/4]$.
\end{enumerate}

\begin{figure}[htp]
\includegraphics[scale=0.40]{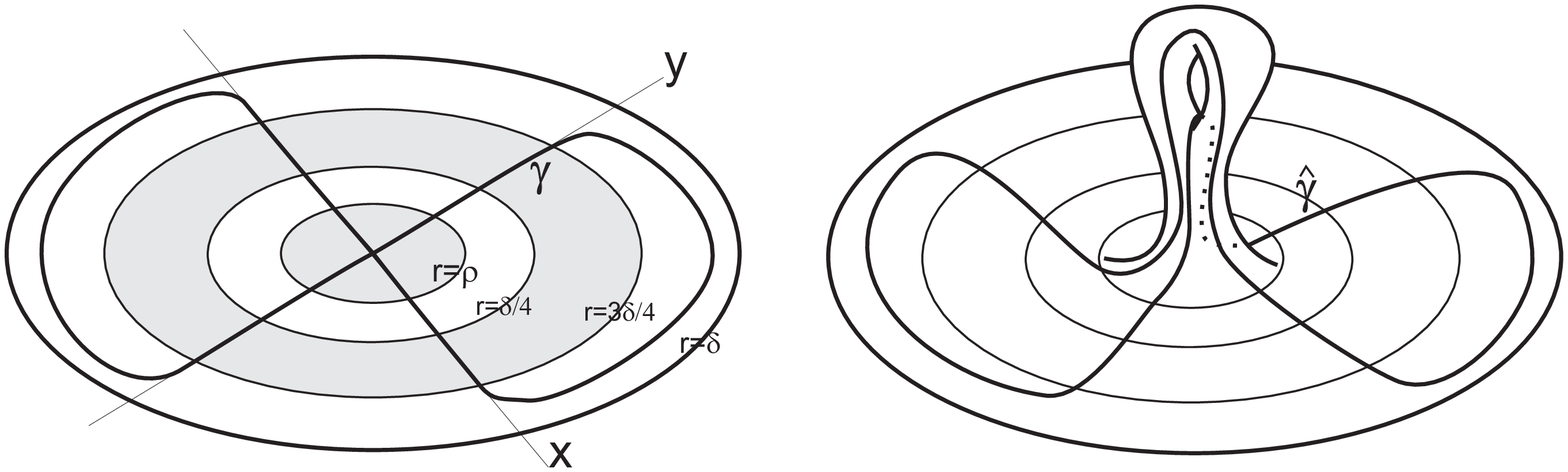} \\
\caption{{\small Gluing the torus in the basis}} \label{fig:glue_tori}
\end{figure}

If we remove a disk, of radius $\rho$, $D^2(\rho) \subset D^2(\delta/4)$ and glue back $T^2 -D^2(\rho)$, we can remove the self-intersection in the path $\gamma$ as shown in the figure \ref{fig:glue_tori}.
\vspace{12mm}
producing a new path $\hat{\gamma}$.

Now we glue the two contact fibrations $M' \supset N \times D^2(\delta) \to D^2(\delta)$ and $N \times T^2 \to T^2$, following the gluing of the basis. They satisfy the conditions to apply Theorem \ref{thm:geiges}. This produces a new contact structure on $M' \bigcup_N (N \times T^2)$. Moreover we obtain a bound in the norms of the contact connection monodromy along $\hat{\gamma}$ in two regions:
\begin{itemize}
\item The gluing area, that is a ring $(-\rho, \rho) \times S^1$, in which we know that a radial path has a monodromy with a generating Hamiltonian $C^2$-bounded by a multiple of $\rho$ because of Lemma \ref{lemma:Geiges_bound}.
\item The torus in which we have a bound given by Lemma \ref{lemma:Frede_bound}. We know that the $C^2$-norm of the Hamiltonian function generating the monodromy is bounded by $C|\epsilon| \delta$.
\end{itemize}

\begin{figure}
\includegraphics[scale=0.40]{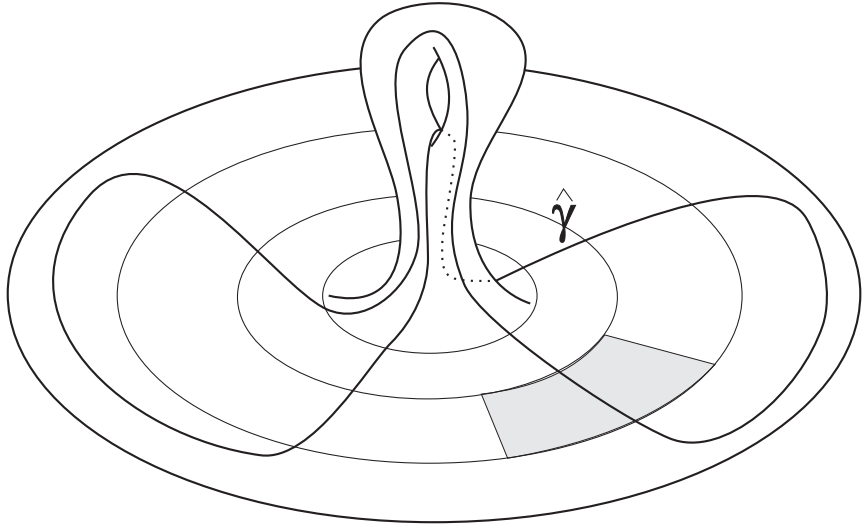}
\caption{{\small Deformation region is depicted in gray.}} \label{fig:pertur_region}
\end{figure}

So we have that $\hat{\gamma}$ produces a parallel transport in the complementary of $\gamma(S^1) \bigcap D^2(\delta)$ with generating Hamiltonian functions ($\hat{\gamma}$ is cut in two disconnected paths in that region) bounded by $C \cdot (\rho + |\epsilon|)$. Therefore we have that the closed path $\gamma$ has a monodromy generated by a Hamiltonian $\{ H_t \}_{t=1}^1$ whose $C^2$-norm is bounded by $C(\rho + |\epsilon|)$. Recall that we are free to choose $\epsilon$ and $\rho$ as small as we wish. Therefore we are in the hypothesis of Proposition \ref{propo:choose_mono} and so a perturbation of the contact fibration on a small region, as shown in Figure \ref{fig:pertur_region}, provides the identity in the monodromy along $\hat{\gamma}$. \hfill $\Box$
\vspace{12mm}

\subsection{Proof of Corollaries \ref{coro:main}, \ref{coro:main2} and \ref{coro:main3}.}
We start  with the $5$-dimensional case. Assume that we want to chek that $S^3 \times \Sigma$ is overtwisted. We take $g(\Sigma)-1$ contact fibrations $S_{ot}^3 \times T^2$ over the torus constructed using Theorem \ref{thm:frederic}. If we fiber-glue them, using Theorem \ref{thm:geiges} we obtain a fibration $M'=S^3_{ot} \times \Sigma'$, where the genus of $\Sigma'$ is $g(\Sigma)-1$. Be aware of the orientations of the normal bundles to perform the gluings: it is needed to choose the constant $\epsilon$ in formula (\ref{eq:frede_e}) positive on one case and negative in the rest. This places $M'$ in the hypothesis of Theorem \ref{thm:main}. So we obtain an overtwisted contact structure in $M' \bigcup_{S^3_{ot}} (S^3_{ot} \times T^2)= S_{ot}^3 \times \Sigma$.

Now, we proceed by induction in the dimension. We just assume that we have been able to find a contact manifold on $M_i= S^3 \times \prod_i \Sigma_i$
with a plastikstufe. We want to construct a plastikstufe in $M_i \times \Sigma$, with $g(\Sigma)\geq 2$. We use Theorems \ref{thm:frederic} and \ref{thm:geiges} to produce a contact structure in $M_i \times \Sigma'$, where $g(\Sigma')= g(\Sigma)-1$. Now we are in the hypothesis of Theorem \ref{thm:main} and we get an overtwisted contact structure in $M_i \times \Sigma$.

To prove Corollary \ref{coro:main3}, we start with the exotic contact structure lying in $S^5$ found in Proposition \ref{propo:exotic_sphere}. There is an overtwisted $S^3 \subset S^5$ whose normal bundle is trivial, so we are in the hypothesis of Theorem \ref{thm:main} and we obtain that $S^5 \bigcup_{S^3} (S^3 \times T^2)=M_o$ admits an overtwisted structure.

Finally Corollary \ref{coro:main2} is a consequence of the connected sum theorem for contact manifolds.

\section{Final remarks.} \label{sect:remarks}
The main results of this paper allow us to show some behaviors, expected to happen, in higher dimensional manifolds. A first example is
\begin{corollary} \label{coro:diversity}
The manifold $S^3\times \Sigma$, for $g(\Sigma) \geq 2$, admits two contact structures in the same homotopy class of hyperplane distributions. The first one being holomorphically fillable, the second one not being fillable at all.
\end{corollary}
{\bf Proof.}
Recall that $S^3 \times \Sigma= \partial (B^4 \times \Sigma)$. Moreover $B^4 \times T^2$ is a complex, and so almost-complex, manifold and admits a $2$-dimensional skeleton defined by the descending disks of the Morse function $f(z_1,z_2,z_3)= |(z_1, z_2)|^2+f_0(z_3)$, where $f_0$ is the standard Morse function in $\Sigma$. Therefore, we are in the hypothesis of Eliashberg's Stein characterization theorem \cite{El90}. Recall that that result produces a Stein structure on $B^4 \times \Sigma$ such that the complex structure is homotopical to the initial one. This implies that $S^3 \times \Sigma$ admits a Stein fillable contact structure in the homotopy class of the hyperplane distribution $D_0 \times \Sigma \subset S^3 \times \Sigma$, where $D_0$ is the standard contact form in $S^3$.

On the other hand the standard structure in $S^3$ lies in the same homotopy class that the standard overtwisted structure $D_{ot}$. Then, Corollary \ref{coro:main} produces a non-fillable contact structure $\hat{\alpha}$ in $S^3 \times \Sigma$. This structure $\ker \hat{\alpha}$ is clearly homotopic to $D_{ot} \times \Sigma$, and therefore homotopic to the Stein-fillable one previously constructed. \hfill $\Box$

Recall that, in dimension $5$ every filling is automatically semipositive, and so Theorem \ref{thm:plas_non} restricts all the cases.

It is possible to copy the very same argument in dimension $7$. We obtain
\begin{corollary}
The manifold $S^3\times \Sigma_1 \times \Sigma_2$ admits two contact structures in the same homotopy class of hyperplane distributions for $g(\Sigma_i)\geq 2$. The first one being holomorphically fillable, the second one being not semipositive-fillable at all.
\end{corollary}
{\bf Proof.}
Recall that $S^3 \times \Sigma_1 \times \Sigma_2= \partial (B^4 \times \Sigma_1 \times \Sigma_2)$. From that point we just follow the argument used in Corollary \ref{coro:diversity} to conclude the result. \hfill $\Box$

\end{document}